# REGULARIZED ESTIMATION OF LARGE COVARIANCE MATRICES

By Peter J. Bickel and Elizaveta Levina[1]

*University of California, Berkeley and University of Michigan*

This paper considers estimating a covariance matrix of $p$ variables from $n$ observations by either banding or tapering the sample covariance matrix, or estimating a banded version of the inverse of the covariance. We show that these estimates are consistent in the operator norm as long as $(\log p)/n \to 0$, and obtain explicit rates. The results are uniform over some fairly natural well-conditioned families of covariance matrices. We also introduce an analogue of the Gaussian white noise model and show that if the population covariance is embeddable in that model and well-conditioned, then the banded approximations produce consistent estimates of the eigenvalues and associated eigenvectors of the covariance matrix. The results can be extended to smooth versions of banding and to non-Gaussian distributions with sufficiently short tails. A resampling approach is proposed for choosing the banding parameter in practice. This approach is illustrated numerically on both simulated and real data.

**1. Introduction.** Estimation of population covariance matrices from samples of multivariate data has always been important for a number of reasons. Principal among these are (1) estimation of principal components and eigenvalues in order to get an interpretable low-dimensional data representation (principal component analysis, or PCA); (2) construction of linear discriminant functions for classification of Gaussian data (linear discriminant analysis, or LDA); (3) establishing independence and conditional independence relations between components using exploratory data analysis and testing; and (4) setting confidence intervals on linear functions of the means of the components. Note that (1) requires estimation of the eigenstructure of the covariance matrix while (2) and (3) require estimation of the inverse.

Received September 2006; revised March 2007.
[1]Supported by NSF Grant DMS-05-05424 and NSA Grant MSPF-04Y-120.

The theory of multivariate analysis for normal variables has been well worked out—see, for example [1]. However, it became apparent that exact expressions were cumbersome, and that multivariate data were rarely Gaussian. The remedy was asymptotic theory for large samples and fixed relatively small dimensions. In recent years, datasets that do not fit into this framework have become very common—the data are very high-dimensional and sample sizes can be very small relative to dimension. Examples include gene expression arrays, fMRI data, spectroscopic imaging, numerical weather forecasting, and many others.

It has long been known that the empirical covariance matrix for samples of size $n$ from a $p$-variate Gaussian distribution, $\mathcal{N}_p(\boldsymbol{\mu}, \Sigma_p)$, has unexpected features if both $p$ and $n$ are large. If $p/n \to c \in (0,1)$ and the covariance matrix $\Sigma_p = I$ (the identity), then the empirical distribution of the eigenvalues of the sample covariance matrix $\hat{\Sigma}_p$ follows the Marĉenko–Pastur law [26], which is supported on $((1-\sqrt{c})^2, (1+\sqrt{c})^2)$. Thus, the larger $p/n$, the more spread out the eigenvalues. Further contributions to the theory of extremal eigenvalues of $\hat{\Sigma}_p$ have been made in [3, 14, 32], among others. In recent years, there have been great developments by Johnstone and his students in the theory of the largest eigenvalues [21, 28] and associated eigenvectors [22]. The implications of these results for inference, other than indicating the weak points of the sample covariance matrix, are not clear.

Regularizing large empirical covariance matrices has already been proposed in some statistical applications—for example, as original motivation for ridge regression [17] and in regularized discriminant analysis [12]. However, only recently has there been an upsurge of both practical and theoretical analyses of such procedures—see [10, 13, 18, 25, 33] among others. These authors study different ways of regularization. Ledoit and Wolf [25] consider Steinian shrinkage toward the identity. Furrer and Bengtsson [13] consider "tapering" the sample covariance matrix, that is, gradually shrinking the off-diagonal elements toward zero. Wu and Pourahmadi [33] use the Cholesky decomposition of the covariance matrix to perform what we shall call "banding the inverse covariance matrix" below, and Huang et al. [18] impose $L_1$ penalties on the Cholesky factor to achieve extra parsimony. Other uses of $L_1$ penalty include applying it directly to the entries of the covariance matrix [2] and applying it to loadings in the context of PCA to achieve sparse representation [34]. Johnstone and Lu [22] consider a different regularization of PCA, which involves moving to a sparse basis and thresholding. Fan, Fan and Lv [10] impose sparsity on the covariance via a factor model.

Implicitly these approaches postulate different notions of sparsity. Wu and Pourahmadi's interest focuses, as does ours, on situations where we can expect that $|i-j|$ large implies near independence or conditional (given the intervening indexes) independence of $X_i$ and $X_j$. At the very least our solutions are appropriate for applications such as climatology and spectroscopy, where there is a natural metric on the index set. Huang et al.'s regularization is more flexible but also depends on the order of variables. In other contexts, notably in finance applications, the sparsity implicit in a factor model of Fan et al. is more suitable. Johnstone and Lu's method relies on the principal eigenvectors being sparse in some basis.

The asymptotic frameworks and convergence results, if at all considered, vary among these studies. Wu and Pourahmadi [33] consider convergence in the sense of single matrix element estimates being close to their population values in probability, with $p_n \to \infty$ at a certain rate determined by the spline smoothers they used. Ledoid and Wolf [25] show convergence of their estimator in "normalized" Frobenius norm $\|A\|_F^2/p$ if $p/n$ is bounded, whereas Furrer and Bengtsoon [13] use the Frobenius norm itself, $\|A\|_F^2 = \text{tr}(AA^T)$, which we shall argue below is too big. Fan, Fan and Lv [10] also show that the Frobenius norm is too big in the factor model context, and employ a different norm based on a sequence of covariance matrices, which is closely related to the entropy loss [20]. Johnstone and Lu [22] show convergence of the first principal component of their estimator when $p/n \to \text{const}$.

We have previously studied [5] the behavior of Fisher's discriminant function for classification as opposed to the so-called "naive Bayes" procedure which is constructed under the assumption of independence of the components. We showed that the latter rule continues to give reasonable results for well-conditioned $\Sigma_p$ as long as $\frac{\log p}{n} \to 0$ while Fisher's rule becomes worthless if $p/n \to \infty$. We also showed that using $k$-diagonal estimators of the covariance achieves asymptotically optimal classification errors if $\Sigma_p$ is Toeplitz and $k_n \to \infty$ at a certain rate. However, the performance of the banded estimators was only evaluated in the context of LDA.

In this paper we show how, by either banding the sample covariance matrix or estimating a banded version of the inverse population covariance matrix, we can obtain estimates which are consistent at various rates in the operator norm as long as $\frac{\log p}{n} \to 0$ and $\Sigma_p$ ranges over some fairly natural families. This implies that maximal and minimal eigenvalues of our estimates and $\Sigma_p$ are close. We introduce the banding approach for the covariance matrix and for the Cholesky factor of the inverse in Section 2. In Section 3 we give the main results: description of classes of covariance matrices for which banding makes sense (Section 3.1), convergence and rates results for the banded covariance estimator (Section 3.2), which we generalize to smooth tapering (Section 3.3) and extend to banding the inverse via its Cholesky factor (Section 3.4). In Section 4 we introduce an analogue of the Gaussian white noise model for covariance matrices and show that if our matrices are embeddable in that model and are well conditioned, then our banded approximations are such that the eigenstructures (individual eigenvalues and associated eigenvectors) of the estimate and population covariance are close. Another approximation result not dependent on existence of the limit model is presented as well. In Section 5 we describe a resampling scheme that can be used to choose the banding parameter $k$ in practice. In Section 6 we give some numerical results, from both simulations and real data. Both theory and simulations indicate that the optimal $k$ depends on $p$, $n$, and the amount of dependence in the underlying model. Section 7 concludes with discussion, and the Appendix contains all the proofs.

**2. The model and two types of regularized covariance estimates.** We assume throughout that we observe $\mathbf{X}_1, \ldots, \mathbf{X}_n$, i.i.d. $p$-variate random variables with mean $\mathbf{0}$ and covariance matrix $\Sigma_p$, and write

$$\mathbf{X}_i = (X_{i1}, \ldots, X_{ip})^T.$$

For now we will assume that the $\mathbf{X}_i$ are multivariate normal, and shall argue separately that it suffices for $X_{1j}^2$ to have subexponential tails for all $j$ (see Extension I after Theorem 1). We want to study the behavior of estimates of $\Sigma_p$ as both $p$ and $n \to \infty$. It is well known that the usual MLE of $\Sigma_p$, the sample covariance matrix,

(1) $$\hat{\Sigma}_p = \frac{1}{n} \sum_{i=1}^{n} (\mathbf{X}_i - \bar{\mathbf{X}})(\mathbf{X}_i - \bar{\mathbf{X}})^T$$

behaves optimally as one might expect if $p$ is fixed, converging to $\Sigma_p$ at rate $n^{-1/2}$. However, as discussed in the Introduction, if $p \to \infty$, $\hat{\Sigma}_p$ can behave very badly unless it is "regularized" in some fashion. Here we consider two methods of regularization.

2.1. *Banding the sample covariance matrix.* For any matrix $M = [m_{ij}]_{p \times p}$, and any $0 \le k < p$, define,

$$B_k(M) = [m_{ij}\mathbf{1}(|i-j| \le k)]$$

and estimate the covariance by $\hat{\Sigma}_{k,p} \equiv \hat{\Sigma}_k = B_k(\hat{\Sigma}_p)$. This kind of regularization is ideal in the situation where the indexes have been arranged in such a way that in $\Sigma_p = [\sigma_{ij}]$

$$|i-j| > k \Rightarrow \sigma_{ij} = 0.$$

This assumption holds, for example, if $\Sigma_p$ is the covariance matrix of $Y_1, \ldots, Y_p$, where $Y_1, Y_2, \ldots$ is a finite inhomogeneous moving average process, $Y_t = \sum_{j=1}^{k} a_{t,t-j} \varepsilon_j$, and $\varepsilon_j$ are i.i.d. mean 0. However, banding an arbitrary covariance matrix does not guarantee positive definiteness—see a generalization to general tapering which avoids this problem in Section 3.3.

2.2. *Banding the inverse.* This method is based on the Cholesky decomposition of the inverse which forms the basis of the estimators proposed by Wu and Pourahmadi [33] and Huang et al. [18]. Suppose we have $\mathbf{X} = (X_1, \ldots, X_p)^T$ defined on a probability space, with probability measure $P$, which is $\mathcal{N}_p(\mathbf{0}, \Sigma_p)$, $\Sigma_p \equiv [\sigma_{ij}]$. Let

$$\hat{X}_j = \sum_{t=1}^{j-1} a_{jt} X_t = \mathbf{Z}_j^T \mathbf{a}_j \tag{2}$$

be the $L_2(P)$ projection of $X_j$ on the linear span of $X_1, \ldots, X_{j-1}$, with $\mathbf{Z}_j = (X_1, \ldots, X_{j-1})^T$ the vector of coordinates up to $j-1$, and $\mathbf{a}_j = (a_{j1}, \ldots, a_{j,j-1})^T$ the coefficients. If $j = 1$, let $\hat{X}_1 = 0$. Each vector $\mathbf{a}_j^T$ can be computed as

$$\mathbf{a}_j = (\operatorname{Var}(\mathbf{Z}_j))^{-1} \operatorname{Cov}(X_j, \mathbf{Z}_j). \tag{3}$$

Let the lower triangular matrix $A$ with zeros on the diagonal contain the coefficients $\mathbf{a}_j$ arranged in rows. Let $\varepsilon_j = X_j - \hat{X}_j$, $d_j^2 = \operatorname{Var}(\varepsilon_j)$ and let $D = \operatorname{diag}(d_1^2, \ldots, d_p^2)$ be a diagonal matrix. The geometry of $L_2(P)$ or standard regression theory imply independence of the residuals. Applying the covariance operator to the identity $\varepsilon = (I - A)X$ gives the modified Cholesky decompositions of $\Sigma_p$ and $\Sigma_p^{-1}$:

$$\begin{aligned} \Sigma_p &= (I - A)^{-1} D [(I - A)^{-1}]^T, \\ \Sigma_p^{-1} &= (I - A)^T D^{-1} (I - A). \end{aligned} \tag{4}$$

Suppose now that $k < p$. It is natural to define an approximation to $\Sigma_p$ by restricting the variables in regression (2) to $\mathbf{Z}_j^{(k)} = (X_{\max(j-k,1)}, \ldots, X_{j-1})$, that is, regressing each $X_j$ on its closest $k$ predecessors only. Replacing $\mathbf{Z}_j$ by $\mathbf{Z}_j^{(k)}$ in (3) gives the new coefficients $\mathbf{a}_j^{(k)}$. Let $A_k$ be the $k$-banded lower triangular matrix containing the new vectors of coefficients $\mathbf{a}_j^{(k)}$, and let $D_k = \operatorname{diag}(d_{j,k}^2)$ be the diagonal matrix containing the corresponding residual variances. Population $k$-banded approximations $\Sigma_{k,p}$ and $\Sigma_{k,p}^{-1}$ are obtained by plugging in $A_k$ and $D_k$ in (4) for $A$ and $D$.

Given a sample $\mathbf{X}_1, \ldots, \mathbf{X}_n$, the natural estimates of $A_k$ and $D_k$ are obtained by performing the operations needed under $\hat{P}$, the empirical distribution, that is, plugging in the ordinary least squares estimates of the coefficients in $A_k$ and the corresponding residual variances in $D_k$. In the general case the variables must be centered first. We will refer to these sample estimates as $\tilde{A}_k = [\tilde{a}_{jt}^{(k)}]$, and $\tilde{D}_k = \operatorname{diag}(\tilde{d}_{j,k}^2)$. Plugging them into (4) for $A$ and $D$ gives the final estimates of $\Sigma_p^{-1}$ and $\Sigma_p$ via the Cholesky decomposition, which we will refer to as $\tilde{\Sigma}_{k,p}^{-1}$ and $\tilde{\Sigma}_{k,p}$, respectively.

Note that since $\tilde{A}_k$ is a $k$-banded lower triangular matrix, $\tilde{\Sigma}_k^{-1}$ is $k$-banded nonnegative definite. Its inverse $\tilde{\Sigma}_k$ is in general not banded, and is different from $\hat{\Sigma}_k$. Similarly, $\tilde{\Sigma}_k^{-1}$ is not the same as $B_k(\hat{\Sigma}^{-1})$, which is in any case ill-defined when $p > n$.

**3. Main results.** All our results can be made uniform on sets of covariance matrices which we define in Section 3.1. Banding the covariance matrix is analyzed in Section 3.2 and generalized to tapering in Section 3.3; results on banding the inverse via the Cholesky decomposition are given in Section 3.4. All the results show convergence of estimators in the matrix $L_2$ norm, $\|M\| \equiv \sup\{\|M\mathbf{x}\| : \|\mathbf{x}\| = 1\} = \lambda_{\max}^{1/2}(M^T M)$, which for symmetric matrices reduces to $\|M\| = \max_i |\lambda_i(M)|$.

3.1. *Classes of covariance matrices.* All our sets will be subsets of the set which we shall refer to as *well-conditioned covariance matrices*, $\Sigma_p$, such that, for all $p$,

$$0 < \varepsilon \leq \lambda_{\min}(\Sigma_p) \leq \lambda_{\max}(\Sigma_p) \leq 1/\varepsilon < \infty.$$

Here, $\lambda_{\max}(\Sigma_p)$, $\lambda_{\min}(\Sigma_p)$ are the maximum and minimum eigenvalues of $\Sigma_p$, and $\varepsilon$ is independent of $p$.

As noted in Bickel and Levina [5], examples of such matrices include covariance matrices of $(U_1, \ldots, U_p)^T$ where $\{U_i, i \geq 1\}$ is a stationary ergodic process with spectral density $f$, $0 < \varepsilon \leq f \leq \frac{1}{\varepsilon}$ and, more generally, of $X_i = U_i + V_i$, $i = 1, \ldots$, where $\{U_i\}$ is a stationary process as above and $\{V_i\}$ is a noise process independent of $\{U_i\}$. This model includes the "spiked model" of Paul [28] since a matrix of bounded rank is Hilbert–Schmidt.

We define the first class of positive definite symmetric well-conditioned matrices $\Sigma \equiv [\sigma_{ij}]$ as follows:

$$\mathcal{U}(\varepsilon_0, \alpha, C) = \bigg\{\Sigma : \max_j \sum_i \{|\sigma_{ij}| : |i - j| > k\} \leq C k^{-\alpha} \text{ for all } k > 0,$$

(5)
$$\text{and } 0 < \varepsilon_0 \leq \lambda_{\min}(\Sigma) \leq \lambda_{\max}(\Sigma) \leq 1/\varepsilon_0\bigg\}.$$

Contained in $\mathcal{U}$ for suitable $\varepsilon_0$, $\alpha$, $C$ is the class

$$\mathcal{L}(\varepsilon_0, m, C) = \{\Sigma : \sigma_{ij} = \sigma(i - j) \text{ (Toeplitz) with spectral density } f_\Sigma,$$

$$0 < \varepsilon_0 \leq \|f_\Sigma\|_\infty \leq \varepsilon_0^{-1}, \|f_\Sigma^{(m)}\|_\infty \leq C\},$$

where $f^{(m)}$ denotes the $m$th derivative of $f$. By Grenander and Szegö [16], if $\Sigma$ is symmetric, Toeplitz, $\Sigma \equiv [\sigma(i - j)]$, with $\sigma(-k) = \sigma(k)$, and $\Sigma$ has an absolutely continuous spectral distribution with Radon–Nikodym derivative $f_\Sigma(u)$, which is continuous on $(-1,1)$, then

$$\|\Sigma\| = \sup_u |f_\Sigma(u)|,$$
$$\|\Sigma^{-1}\| = \left[\inf_u |f_\Sigma(u)|\right]^{-1}. \tag{6}$$

Since $\|f_\Sigma^{(m)}\|_\infty \leq C$ implies that

$$|\sigma(t)| \leq Ct^{-m} \tag{7}$$

which in turn implies $\sum_{t>k} |\sigma(t)| \leq C(m-1)^{-1}k^{-m+1}$, we conclude from 6 and 7 that

$$\mathcal{L}(\varepsilon_0, m, C) \subset \mathcal{U}(\varepsilon_0, m-1, C). \tag{8}$$

A second uniformity class of nonstationary covariance matrices is defined by

$$\mathcal{K}(m, C) = \{\Sigma : \sigma_{ii} \leq Ci^{-m}, \text{ all } i\}.$$

The bound $C$ independent of dimension identifies any limit as being of "trace class" as operator for $m > 1$.

Although $\mathcal{K}$ is not a well-conditioned class,

$$\mathcal{T}(\varepsilon_0, m_1, m_2, C_1, C_2) \equiv \mathcal{L}(\varepsilon_0, m_1, C_1) \oplus \mathcal{K}(m_2, C_2) \subset \mathcal{U}(\varepsilon, \alpha, C), \tag{9}$$

where $\alpha = \min\{m_1-1, m_2/2-1\}$, $C \leq (C_1/(m_1-1) + C_2/(m_2/2-1)$, $\varepsilon^{-1} \leq \varepsilon_0^{-1} + C_2$. To check claim (9), note that

$$\varepsilon_0 \leq \lambda_{\min}(L) \leq \lambda_{\min}(L+K) \leq \lambda_{\max}(L+K)$$
$$\leq \|L\| + \|K\| \leq \varepsilon_0^{-1} + C_2$$

and

$$\max_{j \geq k} \sum_{i:|i-j|>k} |K_{ij}| \leq \max_{j \geq k} \sum_{i:|i-j|>k} |K_{ii}|^{1/2}|K_{jj}|^{1/2}$$
$$\leq C_2(m_2/2-1)^{-1} k^{-m_2/2+1},$$
$$\max_{j<k} \sum_{i:|i-j|>k} |K_{ii}|^{1/2}|K_{jj}|^{1/2} \leq C_2^{1/2} \sum_{i=k+2}^{p} |K_{ii}|^{1/2}$$
$$\leq C_2(m_2/2-1)(k+2)^{-m_2/2+1}.$$

We will use the $\mathcal{T}$ and $\mathcal{L}$ classes for $\Sigma_p$ and $\Sigma_p^{-1}$ for convenience.

3.2. *Banding the covariance matrix.* Our first result establishes rates of convergence for the banded covariance estimator.

THEOREM 1. *Suppose that* **X** *is Gaussian and* $\mathcal{U}(\varepsilon_0, \alpha, C)$ *is the class of covariance matrices defined above. Then, if* $k_n \asymp (n^{-1} \log p)^{-1/(2(\alpha+1))}$,

$$\|\hat{\Sigma}_{k_n,p} - \Sigma_p\| = O_P\left(\left(\frac{\log p}{n}\right)^{\alpha/(2(\alpha+1))}\right) = \|\hat{\Sigma}_{k_n,p}^{-1} - \Sigma_p^{-1}\| \tag{10}$$

*uniformly on* $\Sigma \in \mathcal{U}$.

The proof can be found in the Appendix. Note that the optimal $k_n$ in general depends not only on $p$ and $n$, but also on the dependence structure of the model, expressed by $\alpha$. An approach to choosing $k$ in practice is discussed in Section 5.

From Theorem 1, we immediately obtain:

COROLLARY 1. *If* $\alpha = \min\{m_1 - 1, \frac{m_2}{2} - 1\}$, $m_1 > 1$, $m_2 > 2$, *then* (10) *holds uniformly for* $\Sigma \in \mathcal{T}(\varepsilon_0, m_1, m_2, C_1, C_2)$.

*Extensions of Theorem* 1. I. The Gaussian assumption may be replaced by the following. Suppose $\mathbf{X}_i = (X_{i1}, \ldots, X_{ip})^T$ are i.i.d., $X_{1j} \sim F_j$, where $F_j$ is the c.d.f. of $X_{1j}$, and $G_j(t) = F_j(\sqrt{t}) - F_j(-\sqrt{t})$ is the c.d.f. of $X_{1j}^2$. Then for Theorem 1 to hold it suffices to assume that

$$\max_{1 \le j \le p} \int_0^\infty \exp(\lambda t) \, dG_j(t) < \infty \qquad \text{for } 0 < |\lambda| < \lambda_0 \tag{11}$$

for some $\lambda_0 > 0$. This follows by using the argument of Lemma A.3 and verifying condition (P) on page 45 of [29].

II. If we only assume $E|X_{ij}|^\beta \le C$, $\beta > 2$, for all $j$, we can replace (A4) in the proof of Theorem 1 by

$$P[\|B_k(\hat{\Sigma}^0) - B_k(\Sigma_p)\|_\infty \ge t] \le C n^{-\beta/4} (2k+1) p t^{-\beta/2}. \tag{12}$$

Then a few appropriate modifications of the proof (details omitted here) imply that if $k_n \asymp (n^{-1/2} p^{2/\beta})^{-\gamma(\alpha)}$ where $\gamma(\alpha) = (1 + \alpha + 2/\beta)^{-1}$, then,

$$\|B_{k_n}(\hat{\Sigma}) - \Sigma_p\| = O_P((n^{-1/2} p^{2/\beta})^{\alpha \gamma(\alpha)}). \tag{13}$$

The rate of $k_n$ is still asymptotically optimal.

*Remarks on convergence rates.* (1) Theorem 1 implies that $\|B_{k_n}(\hat{\Sigma}) - \Sigma_p\| \xrightarrow{P} 0$ if $\frac{\log p}{n} \to 0$, uniformly on $\mathcal{U}$. It is not hard to see that if $\Sigma_p = S + K$ where $S$ is Toeplitz, $\varepsilon_0 \le f_S \le \varepsilon_0^{-1}$ and $K$ is trace class in the sense of Section 4, $\Sigma_i K(i,i) < \infty$, then, if $\frac{\log p}{n} \to 0$, there exist $k_n \uparrow \infty$ such that, for the given $\{\Sigma_p\}$,

$$\|B_{k_n}(\hat{\Sigma}) - \Sigma_p\| + \|[B_{k_n}(\hat{\Sigma})]^{-1} - \Sigma_p^{-1}\| \xrightarrow{P} 0. \tag{14}$$

(2) The same claim can be made under (11). On the other hand, under only the moment bound of Extension II with $Ee^{\lambda X_{ij}^2} = \infty$, $\lambda > 0$ we may only conclude that (14) holds if

$$\frac{p^{4/\beta}}{n} \to 0. \tag{15}$$

Related results of Furrer and Bengtsson [13] necessarily have rates of the type (15) not because of tail conditions on the variables, but because they consider the Frobenius norm.

3.3. *General tapering of the covariance matrix.* One problem with simple banding of the covariance matrix is the lack of assured positive definiteness. However, Furrer and Bengtsson [13] have pointed out that positive definiteness can be preserved by "tapering" the covariance matrix, that is, replacing $\hat{\Sigma}_p$ with $\hat{\Sigma}_p * R$, where $*$ denotes Schur (coordinate-wise) matrix multiplication, and $R = [r_{ij}]$ is a positive definite symmetric matrix, since the Schur product of positive definite matrices is also positive definite. This fact was proved by Schur [30] and is also easily seen via a probabilistic interpretation: if $X$ and $Y$ are independent, mean 0 random vectors with $\text{Cov}(\mathbf{X}) = A$, $\text{Cov}(\mathbf{Y}) = B$, then $\text{Cov}(\mathbf{X} * \mathbf{Y}) = A * B$.

In the general case, let $A$ be a countable set of labels of cardinality $|A|$. We can think of a matrix as $[m_{ab}]_{a \in A, b \in A}$. Let $\rho: A \times A \to R^+$, $\rho(a,a) = 0$ for all $a$, be a function we can think of as distance of the point $(a,b)$ from the diagonal. As an example think of $a$ and $b$ as identified with points in $R^m$ and $\rho(a,b) = |a-b|$ where $|\cdot|$ is a norm on $R^m$.

Now suppose $R = [r_{ab}]_{a,b \in A}$ is symmetric positive definite with $r_{ab} = g(\rho(a,b))$, $g: R^+ \to R^+$. Suppose further that $g(0) = 1$ and $g$ is decreasing to 0. Then $R * M$ is a regularization of $M$. Note that $g(t) = 1(t \leq k)$, $\rho(i,j) = |i-j|$ gives banding (which is not nonnegative definite).

In general, let $R_\sigma = [r_\sigma(a,b)]$, where

$$r_\sigma(a,b) = g\left(\frac{\rho(a,b)}{\sigma}\right), \qquad \sigma \geq 0.$$

ASSUMPTION A. *$g$ is continuous, $g(0) = 1$, $g$ is nonincreasing, $g(\infty) = 0$.* Examples of such positive definite symmetric $R_\sigma$ include

$$r_\sigma(i,j) = \left(1 - \frac{|i-j|}{\sigma}\right)_+ \quad \text{and} \quad r_\sigma(i,j) = e^{-|i-j|/\sigma}.$$

With this notation define

$$R_\sigma(M) \equiv [m_{ab} r_\sigma(a,b)]$$

with $R_0(M) = M$. Clearly, as $\sigma \to \infty$, $R_\sigma(M) \to M$.

Our generalization is the following. Denote the range of $g_\sigma(\rho(a,b))$ by $\{g_\sigma(\rho_1), \ldots, g_\sigma(\rho_L)\}$ where $\{0 < \rho_1 < \cdots < \rho_L\}$ is the range of $\rho(a,b)$, $a \in A$, $b \in A$. Note that $L$ depends on $|A| = p$.

THEOREM 2.  *Let $\Delta(\sigma^\varepsilon) = \sum_{l=1}^{L} g_\sigma(\rho_l)$. Note that $\Delta$ depends on $|\mathcal{A}| = p$ and the range of $\rho$. Suppose Assumption A holds. Then if*

$$\Delta \asymp (n^{-1} \log p)^{-1/2(\alpha+1)},$$

*the conclusion of Theorem 1 holds for $R_\sigma(\hat{\Sigma})$.*

The proof of Theorem 2 closely follows the proof of Theorem 1 with (A3) replaced by Lemma A.1 in the Appendix. Both the result and the lemma are of independent interest. The remarks after Theorem 1 generalize equally.

3.4. *Banding the Cholesky factor of the inverse.* Theorems 1 and 2 give the scope of what can be accomplished by banding the sample covariance matrix. Here we show that "banding the inverse" yields very similar results.

If $\Sigma^{-1} = T(\Sigma)^T D^{-1}(\Sigma) T(\Sigma)$ with $T(\Sigma)$ lower triangular, $T(\Sigma) \equiv [t_{ij}(\Sigma)]$, let

$$\mathcal{U}^{-1}(\varepsilon_0, C, \alpha) = \bigg\{\Sigma : 0 < \varepsilon_0 \leq \lambda_{\min}(\Sigma) \leq \lambda_{\max}(\Sigma) \leq \varepsilon_0^{-1},$$

$$\max_i \sum_{j < i-k} |t_{ij}(\Sigma)| \leq C k^{-\alpha} \text{ for all } k \leq p-1 \bigg\}.$$

THEOREM 3. *Uniformly for $\Sigma \in \mathcal{U}^{-1}(\varepsilon_0, C, \alpha)$, if $k_n \asymp (n^{-1} \log p)^{-1/2(\alpha+1)}$ and $n^{-1} \log p = o(1)$,*

$$\|\tilde{\Sigma}_{k_n,p}^{-1} - \Sigma_p^{-1}\| = O_P\bigg(\bigg(\frac{\log p}{n}\bigg)^{\alpha/2(\alpha+1)}\bigg) = \|\tilde{\Sigma}_{k_n,p} - \Sigma_p\|.$$

The proof is given in the Appendix. Note that the condition $n^{-1} \log p = o(1)$ is needed solely for the purpose of omitting a cumbersome and uninformative term from the rate (see Lemma A.2 in the Appendix for details).

It is a priori not clear what $\Sigma \in \mathcal{U}^{-1}$ means in terms of $\Sigma$. The following corollary to Theorem 3 gives a partial answer.

COROLLARY 2. *For $m \geq 2$, uniformly on $\mathcal{L}(\varepsilon_0, m, C)$, if $k_n \asymp (n^{-1} \times \log p)^{-1/2m}$,*

$$\|\tilde{\Sigma}_{k_n,p}^{-1} - \Sigma_p^{-1}\| = O_P\left(\left(\frac{\log p}{n}\right)^{(m-1)/2m}\right) = \|\tilde{\Sigma}_{k_n,p} - \Sigma\|.$$

The proof of Corollary 2 is given in the Appendix. The reason that the argument of Theorem 1 cannot be invoked simply for Theorem 3 is that, as we noted before, $\tilde{\Sigma}^{-1}$ is not the same as $B_k(\hat{\Sigma}^{-1})$, which is not well defined if $p > n$.

**4. An analogue of the Gaussian white noise model and eigenstructure approximations.** In estimation of the means $\boldsymbol{\mu}_p$ of $p$-vectors of i.i.d. variables, the Gaussian white noise model [9] is the appropriate infinite-dimensional model into which all objects of interest are embedded. In estimation of matrices, a natural analogue is the space $\mathcal{B}(l_2, l_2)$, which we write as $\mathcal{B}$, of bounded linear operators from $l_2$ to $l_2$. These can be represented as matrices $[m_{ij}]_{i \geq 1, j \geq 1}$ such that $\Sigma_i [\Sigma_j m_{ij} x_j]^2 < \infty$ for all $\mathbf{x} = (x_1, x_2, \ldots) \in l_2$. It is well known (see Böttcher [7], e.g.) that if $M$ is such an operator, then

$$\|M\|^2 = \sup\{(M\mathbf{x}, M\mathbf{x}) : |\mathbf{x}| = 1\} = \sup \mathcal{S}(M^*M),$$

where $M^*M$ is a self-adjoint member of $\mathcal{B}$ with nonnegative spectrum $\mathcal{S}$. Recall that the spectrum $\mathcal{S}(A)$ of a self-adjoint operator is $\mathcal{R}^c(A)$, where $\mathcal{R}(A) \equiv \{\lambda \in R : A - \lambda J \in \mathcal{B}\}$ where $J$ is the identity. To familiarize ourselves with this space we cite some properties of $\Sigma \in \mathcal{B}$ where

(16) $$\Sigma = [\operatorname{Cov}(X(i), X(j))]_{i,j \geq 1}$$

is the matrix of covariances of a Gaussian stochastic process $\{X(t) : t = 1, 2, \ldots\}$:

1. It is easy to see that the operators $\Sigma$ for all ergodic AR processes, $X(t) = \rho X(t-1) + \varepsilon(t)$ where $\varepsilon(t)$ are i.i.d. $\mathcal{N}(0,1)$ and $|\rho| < 1$, are in $\mathcal{B}$, and $\Sigma^{-1} \in \mathcal{B}$. This is, in fact, true of all ergodic ARMA processes. On the other hand, $X(t) \equiv \sum_{j=1}^{t} \varepsilon(j)$ is evidently not a member of $\mathcal{B}$.

2. The property $\Sigma \in \mathcal{B}$, $\Sigma^{-1} \in \mathcal{B}$ which we shall refer to as being well conditioned, has strong implications. By a theorem of Kolmogorov and Rozanov (see [19]), if $\Sigma$ is Toeplitz, this property holds iff the corresponding stationary Gaussian process is strongly mixing.

We now consider sequences of covariance matrices $\Sigma_p$ such that $\Sigma_p$ is the upper $p \times p$ matrix of the operator $\Sigma \in \mathcal{B}$. That is, $\Sigma$ is the covariance matrix of $\{X(t) : t = 1, 2, \ldots\}$ and $\Sigma_p$ that of $(X(1), \ldots, X(p))$.

By Böttcher [7], if $\Sigma$ is well conditioned, then

$$\Sigma_p(x) \to \Sigma(x)$$

as $p \to \infty$ for all $x \in l_2$. We now combine Theorem 6.1, page 120 and Theorem 5.1, page 474 of Kato [24] to indicate in what ways the spectra and eigenstructures (spectral measures) of $B_{k_n}(\hat{\Sigma}_p)$ are close to those of $\Sigma_p$.

Suppose that the conditions of the remark following Theorem 1 hold. That is, $\Sigma_p$ corresponds to $\Sigma = S + K$, where $S \in \mathcal{B}$ is a Toeplitz operator with spectral density $f_S$ such that, $0 < \varepsilon_0 \leq f_S \leq \varepsilon_0^{-1}$ and $K$ is trace class, $\sum_u K(u,u) < \infty$ which implies $K \in \mathcal{B}$.

Let $M$ be a symmetric matrix and $\mathcal{O}$ be an open set containing $\mathcal{S}(M) \equiv \{\lambda_1, \ldots, \lambda_p\}$ where $\lambda_1(M) \geq \lambda_2(M) \geq \cdots \geq \lambda_p(M)$ are the ordered eigenvalues of $M$ and let $E(M)(\cdot)$ be the spectral measure of $M$ which assigns to each eigenvalue the projection operator corresponding to its eigenspace. Abusing notation, let $E_p \equiv E(\Sigma_p)$, $\hat{E}_p \equiv E(\hat{\Sigma}_{k,p})$, $\mathcal{S} \equiv \mathcal{S}(\Sigma_p)$. Then, $E_p(\mathcal{O}) = E_p(\mathcal{S}) = J$, the identity.

THEOREM 4. *Under the above conditions on $\Sigma_p$,*

$$|\hat{E}_p(\mathcal{O})(x) - x| \xrightarrow{P} 0 \tag{17}$$

*for all $x \in l_2$. Further, if $\mathcal{I}$ is any interval whose endpoints do not belong to $\mathcal{S}$, then*

$$|\hat{E}_p(\mathcal{I} \cap \mathcal{S})(x) - E_p(\mathcal{I})(x)| \xrightarrow{P} 0.$$

Similar remarks apply to $\tilde{\Sigma}_{k,p}$. This result gives no information about rates. It can be refined (Theorem 5.2, page 475 of Kato [24]) but still yields very coarse information. One basic problem is that $\Sigma$ typically has at least in part continuous spectrum and another is that the errors involve the irrelevant bias $|(\Sigma_p - \Sigma)(x)|$. Here is a more appropriate formulation whose consequences for principal component analysis are clear. Let

$$\begin{aligned}\mathcal{G}(\varepsilon, \alpha, C, \Delta, m) \\ = \{\Sigma_p \in \mathcal{U}(\varepsilon, \alpha, C) : \lambda_j(\Sigma_p) - \lambda_{j-1}(\Sigma_p) \geq \Delta, \ 1 \leq j \leq m\}.\end{aligned} \tag{18}$$

Thus the top $m$ eigenvalues are consecutively separated by at least $\Delta$ and all eigenvalues $\lambda_j$ with $j \geq m+1$ are separated from the top $m$ by at least $\Delta$. Furthermore, the dimension of the sum of the eigenspaces of the top $m$ eigenvalues is bounded by $l$ independent of $n$ and $p$. We can then state

THEOREM 5. *Uniformly on $\mathcal{G}$ as above, for $k$ as in Theorem 1, $\mathbf{X}$ Gaussian,*

$$\begin{aligned}|\lambda_j(\hat{\Sigma}_{k,p}) - \lambda_j(\Sigma_p)| &= O_P\bigg(\bigg(\frac{\log p}{n}\bigg)^{1/2}\bigg(\log n + \frac{\alpha}{2}\log p\bigg)\bigg) \\ &= \|E_j(\hat{\Sigma}_{k,p}) - E_j(\Sigma_p)\| \qquad \text{for } 1 \leq j \leq m.\end{aligned} \tag{19}$$

That is, the top $m$ eigenvalues and principal components of $\Sigma_p$, if the eigenvalues are all simple, are well approximated by those of $\hat{\Sigma}_{k,p}$. If we make an additional assumption on $\Sigma_p$,

$$(20) \qquad \frac{\sum_{j=m+1}^{p} \lambda_j(\Sigma_p)}{\sum_{j=1}^{p} \lambda_j(\Sigma_p)} \leq \delta,$$

we can further conclude that the top $m$ principal components of $\hat{\Sigma}_{k,p}$ capture $100(1-\delta)\%$ of the variance of $\mathbf{X}$. To verify (20) we need that

$$(21) \qquad \frac{\operatorname{tr}(\hat{\Sigma}_p - \Sigma_p)}{\operatorname{tr}(\Sigma_p)} = o_P(1).$$

This holds if, for instance, $\operatorname{tr}(\Sigma_p) = \Omega_P(p)$ which is certainly the case for all $\Sigma_p \in \mathcal{T}$. Then, Theorem 5 follows from Theorem 6.1, page 120 of Kato [24], for instance. For simplicity, we give a self-contained proof in the Appendix.

**5. Choice of the banding parameter.** The results in Section 3 give us the rate of $k = k_n$ that guarantees convergence of the banded estimator $\hat{\Sigma}_k$, but they do not offer much practical guidance for selecting $k$ for a given dataset. The standard way to select a tuning parameter is to minimize the risk

$$(22) \qquad R(k) = E\|\hat{\Sigma}_k - \Sigma\|_{(1,1)},$$

with the "oracle" $k$ given by

$$(23) \qquad k_0 = \arg\min_k R(k).$$

The choice of matrix norm in (22) is somewhat arbitrary. In practice, we found the choice of $k$ is not sensitive to the choice of norm; the $l_1$ to $l_1$ matrix norm does just slightly better than others in simulations, and is also faster to compute.

We propose a resampling scheme to estimate the risk and thus $k_0$: divide the original sample into two samples at random and use the sample covariance matrix of one sample as the "target" to choose the best $k$ for the other sample. Let $n_1$, $n_2 = n - n_1$ be the two sample sizes for the random split, and let $\hat{\Sigma}_1^\nu$, $\hat{\Sigma}_2^{(\nu)}$ be the two sample covariance matrices from the $\nu$th split, for $\nu = 1, \ldots, N$. Alternatively, $N$ random splits could be replaced by $K$-fold cross-validation. Then the risk (22) can be estimated by

$$(24) \qquad \hat{R}(k) = \frac{1}{N} \sum_{\nu=1}^{N} \|B_k(\hat{\Sigma}_1^{(\nu)}) - \hat{\Sigma}_2^{(\nu)}\|_{(1,1)}$$

and $k$ is selected as

$$(25) \qquad \hat{k} = \arg\min_k \hat{R}(k).$$

Generally we found little sensitivity to the choice of $n_1$ and $n_2$, and used $n_1 = n/3$ throughout this paper. If $n$ is sufficiently large, another good choice (see, e.g., Bickel, Ritov and Zakai [6]) is $n_1 = n(1 - 1/\log n)$.

The oracle $k_0$ provides the best choice in terms of expected loss, whereas $\hat{k}$ tries to adapt to the data at hand. Another, and more challenging, comparison is that of $\hat{k}$ to the best band choice for the sample in question:

$$(26) \qquad k_1 = \arg\min_k \|\hat{\Sigma}_k - \Sigma\|_{(1,1)}.$$

Here $k_1$ is a random quantity, and its loss is always smaller than that of $k_0$. The results in Section 6 show that $\hat{k}$ generally agrees very well with both $k_0$ and $k_1$, which are quite close for normal data. For heavier-tailed data, one would expect more variability; in that case, the agreement between $\hat{k}$ and $k_1$ is more important that between $\hat{k}$ and $k_0$.

It may be surprising that using the sample covariance $\hat{\Sigma}_2$ as the target in (24) works at all, since it is known to be a very noisy estimate of $\Sigma$. It is, however, an unbiased estimate, and we found that even though (24) tends to overestimate the actual value of the risk, it gives very good results for choosing $k$.

Criterion (24) can be used to select $k$ for the Cholesky-based $\tilde{\Sigma}_k$ as well. An obvious modification—replacing the covariance matrices with their inverses in (24)—avoids additional computational cost and instability associated with computing inverses. One has to keep in mind, however, that while $\hat{\Sigma}_k$ is always well-defined, $\tilde{\Sigma}_k$ is only well-defined for $k < n$, since otherwise regressions become singular. Hence, if $p > n$, $k$ can only be chosen from the range $0, \ldots, n-1$, not $0, \ldots, p-1$.

**6. Numerical results.** In this section, we investigate the performance of the proposed banded estimator of the covariance $\hat{\Sigma}_k$ and the resampling scheme for the choice of $k$, by simulation and on a real dataset. The Cholesky-based $\tilde{\Sigma}_k$ and its variants have been numerically investigated by extensive simulations in [33] and [18], and shown to outperform the sample covariance matrix. Because of that, we omit $\tilde{\Sigma}_k$ from simulations, and only include it in the real data example.

6.1. *Simulations.* We start from investigating the banded estimator by simulating data from $\mathcal{N}(0, \Sigma_p)$ with several different covariance structures $\Sigma_p$. For all simulations, we report results for $n = 100$ and $p = 10$, 100, and 200. Qualitatively, these represent three different cases: $p \ll n$, $p \sim n$ and $p > n$. We have also conducted selected simulations with $p = 1000$, $n = 100$, which qualitatively corresponds to the case $p \gg n$; all the patterns observed with $p > n$ remain the same, only more pronounced. The number of random splits used in (24) was $N = 50$, and the number of replications was 100.

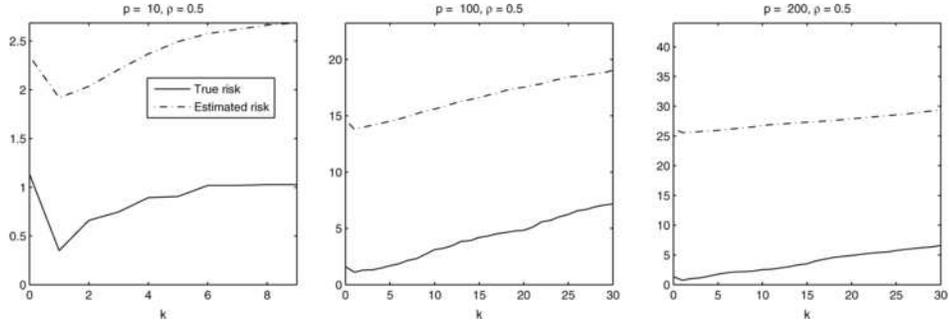

FIG. 1. MA*(1) covariance*: True (*averaged over 100 realizations*) *and estimated risk* (*single realization*) *as a function of $k$, plotted for $k \leq 30$. Both risks are increasing after $k = 1$ for all $p$.*

EXAMPLE 1 (*Moving average covariance structure*). We take $\Sigma_p$ to be the covariance of the MA(1) process, with

$$\sigma_{ij} = \rho^{|i-j|} \cdot \mathbf{1}\{|i-j| \leq 1\}, \qquad 1 \leq i, j \leq p.$$

The true $\Sigma_p$ is banded, and the oracle $k_0 = 1$ for all $p$. For this example we take $\rho = 0.5$. Figure 1 shows plots of the true risk $R(k)$ and the estimated risk $\hat{R}(k)$ from (24). While the risk values themselves are overestimated by (24) due to the extra noise introduced by $\hat{\Sigma}_2$, the agreement of the minima is very good, and that is all that matters for selecting $k$.

Table 1 shows the oracle values of $k_0$ and $k_1$, the estimated $\hat{k}$, and the losses corresponding to all these along with the loss of the sample covariance $\hat{\Sigma}$. When the true model is banded, the estimation procedure always picks the right banding parameter $k = 1$, and performs exactly as well as the oracle. The covariance matrix, as expected, does worse.

EXAMPLE 2 (*Autoregressive covariance structure*). Let $\Sigma_p$ be the covariance of an AR(1) process,

$$\sigma_{ij} = \rho^{|i-j|}, \qquad 1 \leq i, j \leq p.$$

TABLE 1
MA*(1): Oracle and estimated $k$ and the corresponding loss values*

| | | Mean (SD) | | | Loss | | | |
|---|---|---|---|---|---|---|---|---|
| $p$ | $k_0$ | $k_1$ | $\hat{k}$ | $k_1 - \hat{k}$ | $\hat{\Sigma}_{\hat{k}}$ | $\hat{\Sigma}_{k_0}$ | $\hat{\Sigma}_{k_1}$ | $\hat{\Sigma}$ |
| 10 | 1 | 1 (0) | 1 (0) | 0 (0) | 0.5 | 0.5 | 0.5 | 1.2 |
| 100 | 1 | 1 (0) | 1 (0) | 0 (0) | 0.8 | 0.8 | 0.8 | 10.6 |
| 200 | 1 | 1 (0) | 1 (0) | 0 (0) | 0.9 | 0.9 | 0.9 | 20.6 |

For this simulation example, we take $\rho = 0.1$, 0.5 and 0.9. The covariance matrix is not sparse, but the entries decay exponentially as one moves away from the diagonal. Results in Figure 2 and Table 2 show that the smaller $\rho$ is, the smaller the optimal $k$. Results in Table 2 also show the variability in $\hat{k}$ increases when the truth is far from banded (larger $\rho$), which can be expected from the flat risk curves in Figure 2. Variability of $k_1$ increases as well, and $k_1 - \hat{k}$ is not significantly different from 0. In terms of the loss, the estimate again comes very close to the oracle.

EXAMPLE 3 (*Long-range dependence*). This example is designed to challenge the banded estimator, since conditions (5) will not hold for covariance matrix of a process exhibiting long-range dependence. Fractional Gaussian noise (FGN), the increment process of fractional Brownian motion, provides a classic example of such a process. The covariance matrix is given by

$$\sigma_{ij} = \tfrac{1}{2}[(|i-j|+1)^{2H} - 2|i-j|^{2H} + (|i-j|-1)^{2H}], \qquad 1 \leq i,j \leq p,$$

where $H \in [0.5, 1]$ is the Hurst parameter. $H = 0.5$ corresponds to white noise, and the larger $H$, the more long-range dependence in the process. Values of $H$ up to 0.9 are common in practice, for example, in modeling Internet network traffic. For simulating this process, we used the circulant matrix embedding method [4], which is numerically stable for large $p$.

Results in Table 3 show that the procedure based on the estimated risk correctly selects a large $k$ ($k \approx p$) when the covariance matrix contains strong long-range dependence ($H = 0.9$). In this case banding cannot help—but it does not hurt, either, since the selection procedure essentially chooses to do no banding. For smaller $H$, the procedure adapts correctly and selects $k = 0$ for $H = 0.5$ (diagonal estimator for white noise), and a small $k$ for $H = 0.7$.

Another interesting question is how the optimal choice of $k$ depends on dimension $p$. Figure 3 shows the ratio of optimal $k$ to $p$, for both oracle $k_0$ and estimated $\hat{k}$, for AR(1) and FGN [for MA(1), the optimal $k$ is always 1]. The plots confirm the intuition that (a) the optimal amount of regularization depends on $\Sigma$, and the faster the off-diagonal entries decay, the smaller the optimal $k$; and (b) the same model requires relatively more regularization in higher dimensions.

6.2. *Call center data.* Here we apply the banded estimators $\hat{\Sigma}_k$ and $\tilde{\Sigma}_k$ to the call center data used as an example of a large covariance estimation problem by [18], who also provide a detailed description of the data. Briefly, the data consist of call records from a call center of a major U.S. financial institution. Phone calls were recorded from 7:00 am until midnight every day in 2002, and weekends, holidays and days when equipment was malfunctioning have been eliminated, leaving a total of 239 days. On each day, the 17-hour recording period was divided into 10-minute intervals, and the number of calls in each period, $N_{ij}$, was recorded for each of the days $i = 1, \ldots, 239$ and time periods $j = 1, \ldots, 102$. A standard transformation $x_{ij} = (N_{ij} + 1/4)^{1/2}$ was applied to make the data closer to normal.

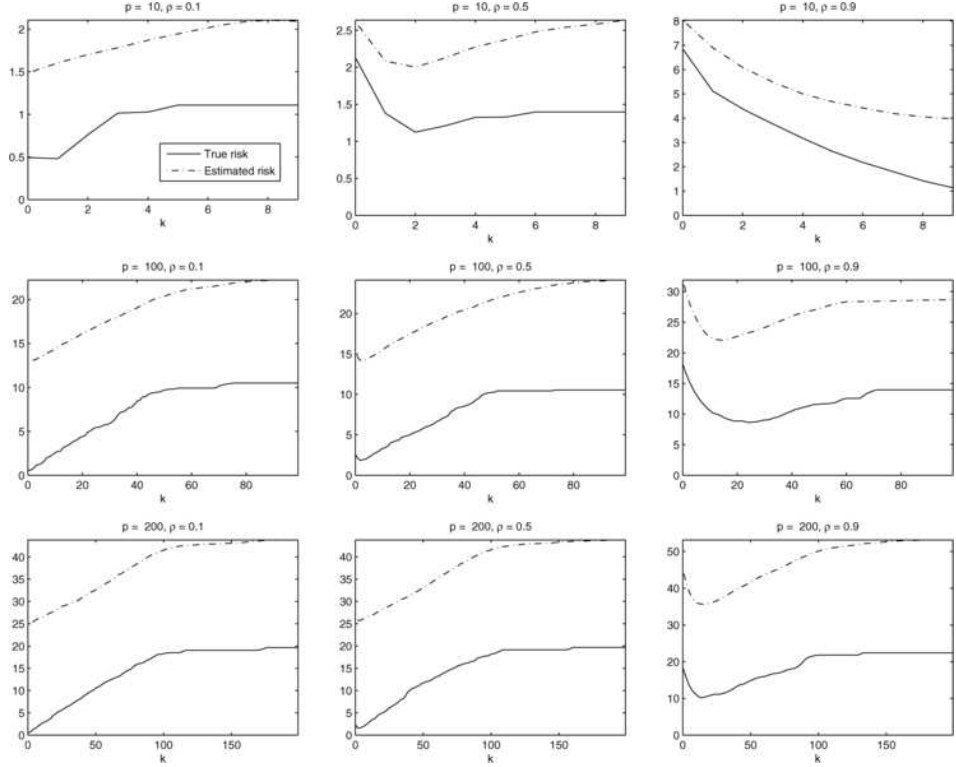

FIG. 2. AR*(1) covariance: True (averaged over 100 realizations) and estimated risk (single realization) as a function of $k$.*

TABLE 2
AR*(1): Oracle and estimated $k$ and the corresponding loss values*

| | | | **Mean (SD)** | | | **Loss** | | | |
|---|---|---|---|---|---|---|---|---|---|
| $p$ | $\rho$ | $k_0$ | $k_1$ | $\hat{k}$ | $k_1 - \hat{k}$ | $\hat{\Sigma}_{\hat{k}}$ | $\hat{\Sigma}_{k_0}$ | $\hat{\Sigma}_{k_1}$ | $\hat{\Sigma}$ |
| 10 | 0.1 | 1 | 0.5 (0.5) | 0.0 (0.2) | 0.5 (0.6) | 0.5 | 0.5 | 0.4 | 1.1 |
| 10 | 0.5 | 3 | 3.3 (0.8) | 2.0 (0.6) | 1.3 (1.1) | 1.1 | 1.0 | 1.0 | 1.3 |
| 10 | 0.9 | 9 | 8.6 (0.7) | 8.9 (0.3) | −0.4 (0.7) | 1.5 | 1.5 | 1.5 | 1.5 |
| 100 | 0.1 | 0 | 0.2 (0.4) | 0.1 (0.3) | 0.1 (0.6) | 0.6 | 0.6 | 0.6 | 10.2 |
| 100 | 0.5 | 3 | 2.7 (0.7) | 2.3 (0.5) | 0.4 (1.0) | 1.6 | 1.6 | 1.5 | 10.6 |
| 100 | 0.9 | 20 | 21.3 (4.5) | 15.9 (2.6) | 5.5 (5.8) | 9.2 | 8.8 | 8.5 | 13.5 |
| 200 | 0.1 | 1 | 0.2 (0.4) | 0.2 (0.4) | −0.0 (0.6) | 0.7 | 0.6 | 0.6 | 20.4 |
| 200 | 0.5 | 3 | 2.4 (0.7) | 2.7 (0.5) | −0.2 (1.0) | 1.8 | 1.7 | 1.7 | 20.8 |
| 200 | 0.9 | 20 | 20.2 (4.5) | 16.6 (2.4) | 3.6 (5.6) | 9.9 | 9.7 | 9.5 | 24.5 |

TABLE 3
*FGN: Oracle and estimated k and the corresponding loss values*

| | | | Mean (SD) | | | $L_1$ Loss | | | |
|---|---|---|---|---|---|---|---|---|---|
| $p$ | $H$ | $k_0$ | $k_1$ | $\hat{k}$ | $k_1 - \hat{k}$ | $\hat{\Sigma}_{\hat{k}}$ | $\hat{\Sigma}_{k_0}$ | $\hat{\Sigma}_{k_1}$ | $\hat{\Sigma}$ |
| 10 | 0.5 | 0 | 0.0 (0.0) | 0.0 (0.0) | 0.0 (0.0) | 0.3 | 0.3 | 0.3 | 1.1 |
| 10 | 0.7 | 5 | 5.0 (1.8) | 2.3 (1.5) | 2.7 (2.5) | 1.4 | 1.2 | 1.1 | 1.2 |
| 10 | 0.9 | 9 | 8.6 (0.6) | 9.0 (0.1) | −0.4 (0.6) | 1.5 | 1.5 | 1.5 | 1.5 |
| 100 | 0.5 | 0 | 0.0 (0.0) | 0.0 (0.0) | 0.0 (0.0) | 0.4 | 0.4 | 0.4 | 10.2 |
| 100 | 0.7 | 4 | 4.9 (2.2) | 4.1 (1.6) | 0.8 (2.9) | 5.5 | 5.5 | 5.4 | 10.7 |
| 100 | 0.9 | 99 | 82.1 (10.9) | 85.1 (15.5) | −3.1 (19.0) | 17.6 | 16.6 | 16.6 | 16.6 |
| 200 | 0.5 | 0 | 0.0 (0.0) | 0.0 (0.1) | −0.0 (0.1) | 0.4 | 0.4 | 0.4 | 20.1 |
| 200 | 0.7 | 3 | 4.2 (2.2) | 4.9 (2.1) | −0.7 (3.4) | 7.9 | 7.7 | 7.7 | 20.9 |
| 200 | 0.9 | 199 | 164.0 (22.7) | 139.7 (38.9) | 24.3 (47.4) | 37.8 | 33.3 | 33.3 | 33.3 |

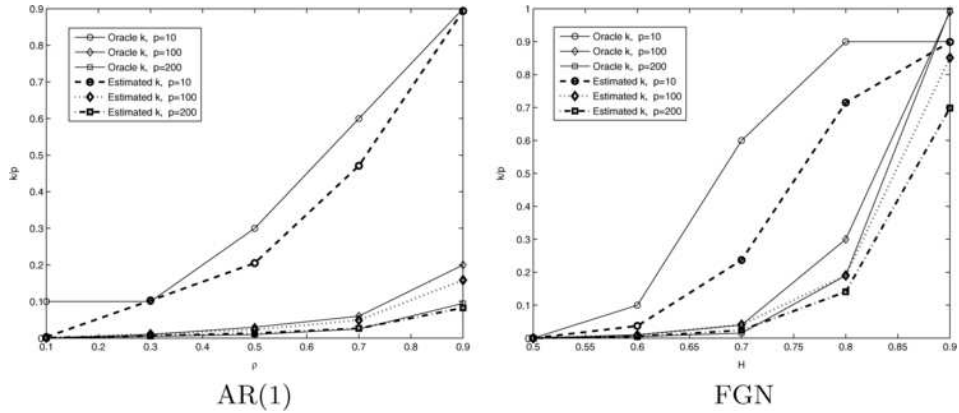

FIG. 3. *The ratio of optimal k to dimension p for* AR*(1) (as a function of $\rho$) and FGN (as a function of H).*

The goal is to predict arrival counts in the second half of the day from counts in the first half of the day. Let $\mathbf{x}_i = (\mathbf{x}_i^{(1)}, \mathbf{x}_i^{(2)})$, with $\mathbf{x}_i^{(1)} = (x_{i1}, \ldots, x_{i,51})$, and $\mathbf{x}_i^{(2)} = (x_{i,52}, \ldots, x_{i,102})$. The mean and the variance of $Vx$ are partitioned accordingly,

$$\mu = \begin{pmatrix} \mu_1 \\ \mu_2 \end{pmatrix}, \qquad \Sigma = \begin{pmatrix} \Sigma_{11} & \Sigma_{12} \\ \Sigma_{21} & \Sigma_{22} \end{pmatrix}. \tag{27}$$

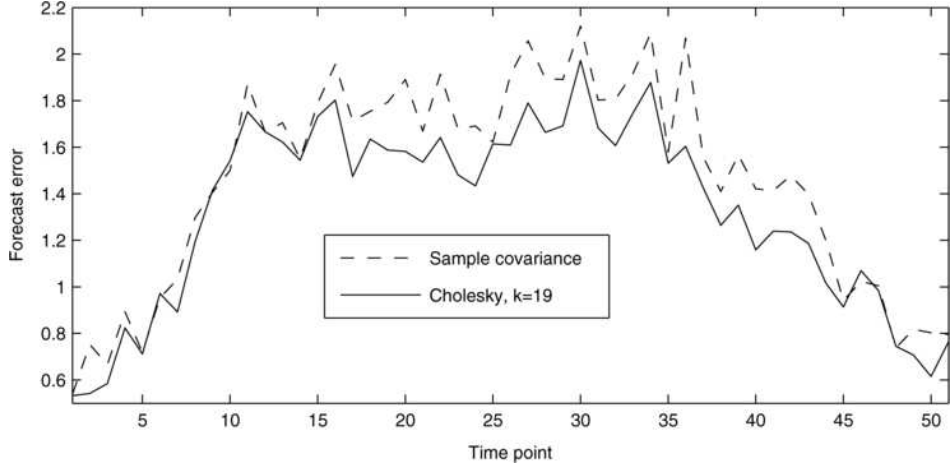

FIG. 4. *Call center forecast error using the sample covariance $\hat{\Sigma}$ and the best Cholesky-based estimator $\tilde{\Sigma}_k$, $k = 19$.*

The best linear predictor of $\mathbf{x}_i^{(2)}$ from $\mathbf{x}_i^{(1)}$ is then given by

$$\hat{\mathbf{x}}_i^{(2)} = \mu_2 + \Sigma_{21}\Sigma_{11}^{-1}(\mathbf{x}_i^{(1)} - \mu_1). \tag{28}$$

Different estimators of $\Sigma$ in (27) can be plugged into (28). To compare their performance, the data were divided into a training set (January to October, 205 days) and a test set (November and December, 34 days). For each time interval $j$, the performance is measured by the average absolute forecast error

$$E_j = \tfrac{1}{34} \sum_{i=206}^{239} |\hat{x}_{ij} - x_{ij}|.$$

The selection procedure for $k$ described in Section 5 to both $\hat{\Sigma}_k$ and $\tilde{\Sigma}_k$. It turns out that the data exhibit strong long-range dependence, and for $\hat{\Sigma}_k$ the selection procedure picks $k = p = 102$, so banding the covariance matrix is not beneficial here. For $\tilde{\Sigma}_k$, the selected $k = 19$ produces a better prediction for almost every time point than the sample covariance $\hat{\Sigma}$ (see Figure 4).

This example suggests that a reasonable strategy for choosing between $\hat{\Sigma}_k$ and $\tilde{\Sigma}_k$ in practice is to estimate the optimal $k$ for both and use the one that selects a smaller $k$. The two estimators are meant to exploit different kinds of sparsity in the data, and a smaller $k$ selected for one of them indicates that that particular kind of sparsity is a better fit to the data.

**7. Discussion.** I. If $\sigma^{ij} = 0$, $|i - j| > k$ and $\|\Sigma^{-1}\| \leq \varepsilon_0^{-1}$, then $\mathbf{X}$ is a $k$th-order autoregressive process and as we might expect, $\tilde{\Sigma}_{k,p}$ is the right

estimate. Now suppose $\sigma^{ii} \leq \varepsilon_0^{-1}$ for all $i$ and we only know that $\sigma^{ij} = 0$ for each $i$ and $p - (2k+1)$ $j$'s. This condition may be interpreted as saying that, for each $i$ there is a set $S_i$ with $|S_i| \leq k$, $i \notin S_0$, such that, $X_i$ is independent of $\{X_t, t \notin S_i, t \neq i\}$ given $\{X_j : j \in S_i\}$. Although banding would not in general give us sparse estimates, the following seem intuitively plausible:

(1) Minimize a suitable objective function $\Psi(\hat{P}, \Sigma) \geq 0$ where $\hat{P}$ is the empirical distribution of $\mathbf{X}_1, \ldots, \mathbf{X}_n$ and

$$\Psi(P, \Sigma_p) = 0$$

subject to $\|\Sigma\|_{(1,1)} \leq \gamma_{n,p}$.

(2) Let $\gamma_{n,p} \to 0$ "slowly." This approach should yield estimates which consistently estimate sparse covariance structure. Banerjee, D'Aspremont and El Ghaoui [2] and Huang et al. [18] both use normal or Wishart-based log-likelihoods for $\Psi$ and a Lasso-type penalty in this context. We are currently pursuing this approach more systematically.

II. The connections with graphical models are also apparent. If $D$ is the dependency matrix of $\Sigma^{-1}$, with entries 0 and 1, then $\|D\|_{(1,1)}$ is just the maximum degree of the graph vertices. See Meinshausen and Buhlmann [27] for a related approach in determining covariance structure in this context.

III. A similar interpretation can be attached if we assume $\Sigma$ is $k_0$ banded after a permutation of the rows. This is equivalent to assuming that there is a permutation of variables after which $X_i$ is independent of $\{X_j : |j - i| > k\}$ for all $i$.

IV. In the case when $\Sigma$ is the covariance of a stationary process, there is a possibility of regularizing the covariance function via inverting a regularized estimate of the spectral density. There is a large literature on spectral density estimate (see, e.g., [11] for an approach based on local smoothing and a review). Exploring this connection further, and in particular, understanding the equivalent of banding in the spectral domain, is a topic for future work.

## APPENDIX: ADDITIONAL LEMMAS AND PROOFS

In addition to the operator norm $\|M\|$ from $l_2$ to $l_2$ we have defined, in what follows we use additional matrix norms. For a vector $\mathbf{x} = (x_1, \ldots, x_p)^T$, let

$$\|\mathbf{x}\|_1 = \sum_{j=1}^{p} |x_j|, \qquad \|\mathbf{x}\|_\infty = \max_j |x_j|.$$

For a matrix $M = [m_{ij}]$, the corresponding operator norms from $l_1$ to $l_1$ and $l_\infty$ to $l_\infty$ are, respectively,

$$\|M\|_{(1,1)} \equiv \sup\{\|M\mathbf{x}\|_1 : \|\mathbf{x}\|_1 = 1\} = \max_j \sum_i |m_{ij}|,$$

(A1)
$$\|M\|_{(\infty,\infty)} \equiv \sup\{\|M\mathbf{x}\|_\infty : \|\mathbf{x}\|_\infty = 1\} = \max_i \sum_j |m_{ij}|.$$

We will also write $\|M\|_\infty \equiv \max_{i,j} |m_{ij}|$.

For symmetric matrices, $\|M\|_{(1,1)} = \|M\|_{(\infty,\infty)}$. The $l_1$ to $l_1$ norm arises naturally through the inequality (see, e.g., Golub and Van Loan [15])

(A2) $\|M\| \leq [\|M\|_{(1,1)} \|M\|_{(\infty,\infty)}]^{1/2} = \|M\|_{(1,1)}$ for $M$ symmetric.

PROOF OF THEOREM 1. It is easy to see that (A2) and (A2) imply

(A3) $\qquad \|B_k(\hat{\Sigma}) - B_k(\Sigma)\| = O_P(k\|B_k(\hat{\Sigma}) - B_k(\Sigma)\|_\infty).$

Let $\hat{\Sigma}^0 = \frac{1}{n}\sum_{i=1}^n \mathbf{X}_i^T \mathbf{X}_i$ and w.l.o.g. $E\mathbf{X}_1 = \mathbf{0}$. By an application of a result of Saulis and Statulevičius [29] (see Lemma A.3) and the union sum inequality,

(A4) $\qquad P[\|B_k(\hat{\Sigma}^0) - B_k(\Sigma)\|_\infty \geq t] \leq (2k+1)p \ \exp\{-nt^2\gamma(\varepsilon_0,\lambda)\}$

for $|t| \leq \lambda \equiv \lambda(\varepsilon_0)$. By choosing $t = M(\frac{\log(pk)}{n})^{1/2}$ for $M$ arbitrary we conclude that, uniformly on $\mathcal{U}$,

(A5)
$$\|B_k(\hat{\Sigma}^0) - B_k(\Sigma_p)\|_\infty = O_P((n^{-1}\log(pk))^{1/2})$$
$$= O_P((n^{-1}\log p)^{1/2})$$

since $k < p$. On the other hand, by 5,

(A6) $\qquad \|B_k(\Sigma_p) - \Sigma_p\|_\infty \leq Ck^{-\alpha}$

for $\Sigma_p \in \mathcal{U}(\varepsilon_0, \alpha, C)$.

Combining (A5) and (A6), the result follows for $B_k(\hat{\Sigma}^0)$. But, if $\bar{\mathbf{X}} = (\bar{X}_1, \ldots, \bar{X}_p)^T$,

$$\|B_k(\hat{\Sigma}^0) - B_k(\hat{\Sigma})\| \leq \|B_k(\bar{\mathbf{X}}^T\bar{\mathbf{X}})\| \leq (2k+1)\max_{1\leq j\leq p}|\bar{X}_j|^2$$
$$= O_P\left(\frac{k\log p}{n}\right) = O_P((n^{-1}\log p)^{\alpha/(2(\alpha+1))}).$$

Since
$$\|[B_{k_n}(\hat{\Sigma})]^{-1} - \Sigma_p^{-1}\| = \Omega_P(\|B_{k_n}(\hat{\Sigma}) - \Sigma_p\|),$$

uniformly on $\mathcal{U}$, the result follows. □

The key to Theorem 2 is the following lemma which substitutes for (A3). Consider symmetric matrices $M$ indexed by $(a,b)$, $a,b \in \mathcal{A}$, a finite index set. Suppose for each $a \in \mathcal{A}$ there exist $N_a \leq N$ sets $S_{a,j}$ such that the $S_{a,j}$ form a partition of $\mathcal{A} - \{a\}$. Define, for any $1 \leq j \leq N$, $M = [m(a,b)]$ as above:

$$r(j) = \max\{|m(a,b)|: b \in S_{a,j}, \ a \in \mathcal{A}\}$$

and $\mu = \max_a |m(a,a)|$.

LEMMA A.1. *Under Assumption* A,

$$\|M\| \leq \mu + \sum_{j=1}^{N} r(j). \tag{A7}$$

PROOF. Apply (A2) noting that

$$\sum \{|m(a,b)|: b \in \mathcal{A}\} \leq \sum_{j=1}^{N} r(j) + \mu$$

for all $a \in \mathcal{A}$. □

PROOF OF COROLLARY 2. An examination of the proof of Theorem 1 will show that the bound of $\|\Sigma_p - B_k(\Sigma_p)\|_{(1,1)}$ was used solely to bound $\|\Sigma_p - B_k(\Sigma_p)\|$. But in the case of Corollary 2, a theorem of Kolmogorov (De Vore and Lorentz [8], page 334) has, after the identification (A1),

$$\|\Sigma_p - B_k(\Sigma_p)\| \leq \frac{C' \log k}{k^m} \tag{A8}$$

where $C'$ depends on $C$ and $m$ only, for all $\Sigma_p \in \mathcal{L}(\varepsilon_0, m, C)$. The result follows. Note that Corollary 1 would give the same results as the inferior bound $C'k^{-(m-1)}$. □

To prove Theorem 3 we will need:

LEMMA A.2. *Under conditions of Theorem* 3, *uniformly on* $\mathcal{U}$,

$$\max\{\|\tilde{\mathbf{a}}_j^{(k)} - \mathbf{a}_j^{(k)}\|_\infty: 1 \leq j \leq p\} = O_P(n^{-1/2} \log^{1/2} p), \tag{A9}$$

$$\max\{|\tilde{d}_{j,k}^2 - d_{j,k}^2|: 1 \leq j \leq p\} = O_P((n^{-1} \log p)^{\alpha/(2(\alpha+1))}) \tag{A10}$$

*and*

$$\|A_k\| = \|D_k^{-1}\| = O(1), \tag{A11}$$

*where* $\tilde{\mathbf{a}}_j^{(k)} = (\tilde{a}_{j1}^{(k)}, \ldots, \tilde{a}_{j,j-1}^{(k)})$ *are the empirical estimates of the vectors* $\mathbf{a}_j^{(k)} = (a_{j1}^{(k)}, \ldots, a_{j,j-1}^{(k)})$ *and* $\tilde{d}_{j,k}^2$, $1 \leq j \leq p$ *are the empirical estimates of the* $d_{j,k}^2$.

To prove Lemma A.2 we need an additional lemma, obtained from results of Saulis and Statulevičius [29].

LEMMA A.3. *Let $Z_i$ be i.i.d. $\mathcal{N}(\mathbf{0}, \Sigma_p)$ and $\lambda_{\max}(\Sigma_p) \leq \varepsilon_0^{-1} < \infty$. Then, if $\Sigma_p = [\sigma_{ab}]$,*

$$P\left[\left|\sum_{i=1}^n (Z_{ij}Z_{ik} - \sigma_{jk})\right| \geq n\nu\right]$$

(A12)
$$\leq C_1 \exp(-C_2 n\nu^2) \qquad \text{for } |\nu| \leq \delta,$$

*where $C_1$, $C_2$ and $\delta$ depend on $\varepsilon_0$ only.*

PROOF. Write

$$P\left[\left|\sum_{i=1}^n (Z_{ij}Z_{ik} - \sigma_{jk})\right| \geq n\nu\right]$$

$$= P\left[\left|\sum_{i=1}^n (Z_{ij}^* Z_{ik}^* - \rho_{jk})\right| \geq \frac{n\nu}{(\sigma_{jj}\sigma_{kk})^{1/2}}\right],$$

where $\rho_{jk} = \sigma_{jk}(\sigma_{jj}\sigma_{kk})^{-1/2}$ and $(Z_{ij}^*, Z_{ik}^*) \sim \mathcal{N}_2(0,0,1,1,\rho_{jk})$. Now,

$$\sum_{i=1}^n (Z_{ij}^* Z_{ik}^* - \rho_{jk})$$

(A13)
$$= \tfrac{1}{4}\left[\sum_{i=1}^n [(Z_{ij}^* + Z_{ik}^*)^2 - 2(1+\rho_{jk})]\right.$$
$$\left. + \sum_{i=1}^n [(Z_{ij}^* - Z_{ik}^*)^2 - 2(1-\rho_{jk})]\right]$$

and reduce the problem to estimating

$$2P\left[\left|\sum_{i=1}^n (V_i^2 - 1)\right| \geq \frac{n\nu}{2(1-\rho_{jk})(\sigma_{jj}\sigma_{kk})^{1/2}}\right],$$

where $V_i$ are i.i.d. $\mathcal{N}(0,1)$. Since $\chi_1^2$ satisfies condition (P) (3.12) on page 45 of [29], the lemma follows from Theorem 3.2, page 45 and (2.13) on page 19, since $(\sigma_{jj}\sigma_{kk})^{1/2}|1-\rho_{jk}| \leq 2\varepsilon_0^{-1}$. □

PROOF OF LEMMA A.2. Note first that

(A14) $\qquad \|\operatorname{Var}\mathbf{X} - \widehat{\operatorname{Var}}\mathbf{X}\|_\infty = O_P(n^{-1/2}\log^{1/2} p),$

by Lemma A.3. Hence,

(A15) $\quad \max_j \|\widehat{\operatorname{Var}}^{-1}(\mathbf{Z}_j^{(k)}) - \operatorname{Var}^{-1}(\mathbf{Z}_j^{(k)})\|_\infty = O_P(n^{-1/2}\log^{1/2} p).$

To see this, note that the entries of $\widehat{\operatorname{Var}}\mathbf{X} - \operatorname{Var}\mathbf{X}$ can be bounded by $n^{-1}|\sum_{i=1}^n X_{ia}X_{ib} - \sigma_{ab}| + n^{-2}|\sum_{i=1}^n X_{ia}||\sum_{i=1}^n X_{ib}|$, where w.l.o.g. we assume $E\mathbf{X} = \mathbf{0}$. Lemma A.3 ensures that

$$P\left[\max_{a,b}\left|n^{-1}\sum_{i=1}^n (X_{ia}X_{ib} - \sigma_{ab})\right| \geq \nu\right] \leq C_1 p^2 \exp(-C_2 n \nu^2)$$

for $|\nu| \leq \delta$. Now take $\nu = (\frac{\log p^2}{nC_2})^{1/2} M$ for $M$ arbitrary. The second term is similarly bounded.

Also,

$$\|\Sigma^{-1}\| = \|(\operatorname{Var} X)^{-1}\| \leq \varepsilon_0^{-1}.$$

Claim (A9) and the first part of (A11) follow from 3, (A14) and (A15). Since

$$\tilde{d}_{jk}^2 = \widehat{\operatorname{Var}} X_j - \widehat{\operatorname{Var}}\left(\sum_{t=j-k}^{j-1} \tilde{a}_{jt}^{(k)} X_t\right),$$

$$d_{jk}^2 = \operatorname{Var} X_j - \operatorname{Var}\left(\sum_{t=j-k}^{j-1} a_{jt}^{(k)} X_t\right),$$

and the covariance operator is linear,

$$\begin{aligned}|\tilde{d}_{jk}^2 - d_{jk}^2| &\leq |\operatorname{Var}(X_j) - \widehat{\operatorname{Var}} X_j| \\
&\quad + \left|\widehat{\operatorname{Var}}\sum_{t=j-k}^{j-1}(\tilde{a}_{jt}^{(k)} - a_{jt}^{(k)})X_t\right| \\
&\quad + \left|\widehat{\operatorname{Var}}\sum_{t=j-k}^{j-1} a_{jt}^{(k)} X_t - \operatorname{Var}\sum_{t=j-k}^{j-1} a_{jt}^{(k)} X_t\right|.\end{aligned} \quad \text{(A16)}$$

The sum $\sum_{t=j-k}^{j-1}$ is understood to be $\sum_{t=\max(1,j-k)}^{j-1}$. The maximum over $j$ of the first term is $O_P(n^{-1/2}\log^{1/2} p)$ by Lemma A.3. The second can be written as

$$\begin{aligned}&\left|\sum\{(\tilde{a}_{js}^{(k)} - a_{js}^{(k)})(\tilde{a}_{jt}^{(k)} - a_{jt}^{(k)})\widehat{\operatorname{Cov}}(X_s, X_t) : j - k \leq s,\ t \leq j - 1\}\right| \\
&\leq \left(\sum_{t=j-k}^{j-1}|\tilde{a}_{jt}^{(k)} - a_{jt}^{(k)}|\widehat{\operatorname{Var}}^{1/2}(X_t)\right)^2 \\
&\leq k^2 \max_t (\tilde{a}_{jt}^{(k)} - a_{jt}^{(k)})^2 \max_t \widehat{\operatorname{Var}}(X_t) \\
&= O_P(k^2 n^{-1}(\log p)^2) = O_P((n^{-1}\log p)^{\alpha/(2(\alpha+1))})\end{aligned} \quad \text{(A17)}$$

by (A9) and $\|\Sigma_p\| \leq \varepsilon_0^{-1}$. Note that in the last equality we used the assumption $n^{-1} \log p = o(1)$. The third term in (A16) is bounded similarly. Thus (A10) follows. Further, for $1 \leq j \leq p$,

$$
\begin{aligned}
d_{jk}^2 &= \operatorname{Var}\left(X_j - \sum \{a_{jt}^{(k)} X_t \colon \max(1, j-k) \leq t \leq j-1\}\right) \\
&\geq \varepsilon_0 \left(1 + \sum (a_{jt}^{(k)})^2\right) \geq \varepsilon_0
\end{aligned}
\tag{A18}
$$

and the lemma follows. □

PROOF OF THEOREM 3. We parallel the proof of Theorem 1. We need only check that

$$\|\tilde{\Sigma}_{k,p}^{-1} - \Sigma_{k,p}^{-1}\|_\infty = O_P(n^{-1/2} \log^{1/2} p) \tag{A19}$$

and

$$\|\Sigma_{k,p}^{-1} - B_k(\Sigma_p^{-1})\| = O(k^{-\alpha}). \tag{A20}$$

We first prove (A19). By definition,

$$
\begin{aligned}
&\tilde{\Sigma}_{k,p}^{-1} - \Sigma_{k,p}^{-1} \\
&= (I - \tilde{A}_k)\tilde{D}_k^{-1}(I - \tilde{A}_k)^T - (I - A_k)D_k^{-1}(I - A_k)^T
\end{aligned}
\tag{A21}
$$

where $\tilde{A}_k, \tilde{D}_k$ are the empirical versions of $A_k$ and $D_k$. Apply the standard inequality

$$
\begin{aligned}
&\|A^{(1)} A^{(2)} A^{(3)} - B^{(1)} B^{(2)} B^{(3)}\| \\
&\leq \sum_{j=1}^3 \|A^{(j)} - B^{(j)}\| \prod_{k \neq j} \|B^{(k)}\| \\
&\quad + \sum_{j=1}^3 \|B^{(j)}\| \prod_{k \neq j} \|A^{(k)} - B^{(k)}\| + \prod_{j=1}^3 \|A^{(j)} - B^{(j)}\|.
\end{aligned}
\tag{A22}
$$

Take $A^{(1)} = [A^{(3)}]^T = I - \tilde{A}_k$, $B^{(1)} = [B^{(3)}]^T = I - A_k$, $A^{(2)} = \tilde{D}_k^{-1}$, $B^{(2)} = D_k^{-1}$ in (A22) and (A19) follows from Lemma A.2. For (A20), we need only note that for any matrix $M$,

$$
\begin{aligned}
&\|MM^T - B_k(M)B_k(M^T)\| \\
&\qquad \leq 2\|M\|\|B_k(M) - M\| + \|B_k(M) - M\|^2
\end{aligned}
$$

and (A20) and the theorem follows from our definition of $U^{-1}$. □

LEMMA A.4. *Suppose $\Sigma = [\rho(j-i)]$ is a Toeplitz covariance matrix; $\rho(k) = \rho(-k)$ for all $k$, $\Sigma \in \mathcal{L}(\varepsilon_0, m, C)$. Then, if $f$ is the spectral density of $\Sigma$:*

(i) $\Sigma^{-1} = [\tilde{\rho}(j-i)]$, $\tilde{\rho}(k) = \tilde{\rho}(-k)$,
(ii) $\Sigma^{-1}$ *has spectral density* $\frac{1}{f}$,
(iii) $\Sigma^{-1} \in \mathcal{L}(\varepsilon_0, m, C'(m, \varepsilon_0, C))$.

PROOF. That $\|(\frac{1}{f})^{(m)}\|_\infty \leq C'(m, \varepsilon_0, C)$ and $\varepsilon_0 \leq \|\frac{1}{f}\|_\infty \leq \varepsilon_0^{-1}$ is immediate. The claims (i) and (ii) follow from the identity, $\frac{1}{f} = \sum_{k=-\infty}^{\infty} \tilde{\rho}(k) e^{2\pi i k u}$ in the $L_2$ sense and

$$1 = \sum_{k=-\infty}^{\infty} \delta_{0k} e^{2\pi i k u} = f(u) \frac{1}{f}(u). \qquad \square$$

PROOF OF COROLLARY 2. Note that $\Sigma \in \mathcal{L}(\varepsilon_0, m, C_0)$ implies that

$$(A23) \qquad f_\Sigma^{-1/2}(u) = a_0 + \sum_{j=1}^{\infty} a_k \cos(2\pi j u)$$

is itself $m$ times differentiable and

$$(A24) \qquad \|(f_\Sigma^{-1/2})^{(m)}\|_\infty \leq C'(C_0, \varepsilon_0).$$

But then,

$$(A25) \qquad \begin{aligned} f_{\Sigma^{-1}}(u) &= b_0 + \sum_{j=1}^{\infty} b_j \cos 2\pi j u \\ &= \left(a_0 + \sum_{j=1}^{\infty} a_j \cos 2\pi j u\right)^2, \end{aligned}$$

where $b_i = \sum_{j=0}^{i} a_j a_{i-j}$. All formal operations are justified since $\sum_{j=0}^{\infty} |a_j| < \infty$ follows from Zygmund [35], page 138. But (A25) can be reinterpreted in view of Lemma A.4 as $\Sigma^{-1} = AA^T$ where $A = [a_{i-j} 1(i \geq j)]$ and $a_j$ are real and given by (A25). Then, if $A_k \equiv B_k(A)$, $B_k(A) B_k^T(A)$ has spectral density,

$$(A26) \qquad f_{\Sigma_{k,p}^{-1}}(u) = \left(\sum_{j=0}^{k} a_j \cos 2\pi j u\right)^2.$$

Moreover, from (A25) and (A26)

$$\|f_{\Sigma_{k,p}^{-1}} - f_{\Sigma_p^{-1}}\|_\infty$$

$$\leq \left\|\sum_{j=k+1}^{\infty} a_j \cos 2\pi j u\right\|_\infty \left(\|f_{\Sigma_p}^{-1/2}\|_\infty + \left\|\sum_{j=k+1}^{\infty} a_j \cos 2\pi j u\right\|_\infty\right).$$

By (A24) $|a_j| \leq C' j^{-m}$, hence finally,

$$(A27) \qquad \|\Sigma_{k,p}^{-1} - \Sigma_p^{-1}\| = \|f_{\Sigma_{k,p}^{-1}} - f_{\Sigma^{-1}}\|_\infty \leq C k^{-(m-1)}.$$

Corollary 2 now follows from (A27) and (A9) and (A10) by minimizing

$$C_1 \frac{k^3 \log^{1/2}(pk)}{n^{1/2}} + C_2 k^{-(m-1)}. \qquad \Box$$

PROOF OF THEOREM 5. We employ a famous formula of Kato [23] and Sz.-Nagy [31]. If $R(\lambda, M) \equiv (M - \lambda J)^{-1}$ for $\lambda \in \mathcal{S}^c$, the resolvent set of $M$ and $\lambda_0$ is an isolated eigenvalue, $|\lambda - \lambda_0| \geq \Delta$ for all $\lambda \in \mathcal{S}$, $\lambda \neq \lambda_0$, then (formula (1.16), page 67, Kato [24])

$$(A28) \qquad E_0(x) = \frac{1}{2\pi i} \int_\Gamma R(\lambda, M) \, d\lambda,$$

where $E_0$ is the projection operator on the eigenspace corresponding to $\lambda_0$ and $\Gamma$ is a closed simple contour in the complex plane about $\lambda_0$ containing no other member of $\mathcal{S}$. The formula is valid not just for symmetric $M$ but we only employ it there. We argue by induction on $m$. For $m = 1$, $|\lambda_1(M) - \lambda_1(N)| \leq \|M - N\|$ for $M$, $N$ symmetric by the Courant–Fischer theorem. Thus, if $\|\hat{\Sigma}_{k,p} - \Sigma_p\| \leq \frac{\Delta}{2}$ (say), we can find $\Gamma$ containing $\lambda_1(\hat{\Sigma}_{k,p})$ and $\lambda_1(\Sigma_p)$ and no other eigenvalues of either matrix with all points on $\Gamma$ at distance at least $\Delta/4$ from both $\lambda_1(\hat{\Sigma}_{k,p})$ and $\lambda_1(\Sigma_p)$. Applying (A28) we conclude that

$$\|\hat{E}_1 - E_1\| \leq \max_\Gamma \{\|R(\lambda, \Sigma_p)\| \|R(\lambda, \hat{\Sigma}_{k,p})\|\} \|\hat{\Sigma}_{k,p} - \Sigma\|.$$

But, $\|R(\lambda, \Sigma_p)\| \leq \{\max_j |\lambda - \lambda_j(\Sigma_p)|\}^{-1} \leq 2/\Delta$ by hypothesis, and similarly for $\|R(\lambda, \hat{\Sigma}_{k,p})\|$. Therefore,

$$(A29) \qquad \|\hat{E}_1 - E_1\| \leq 16\Delta^{-2} \|\hat{\Sigma}_{k,p} - \Sigma\|,$$

and the claims (19) and (20) are established for $m = 1$. We describe the induction step from $m = 1$ to $m = 2$ which is repeated with slightly more cumbersome notation for all $m$ (omitted). Consider a unit vector,

$$(A30) \qquad x = \sum_{j=2}^p E_j x \perp E_1 x = (\hat{E}_1 - E_1)x + (J - \hat{E}_1)x.$$

Then,

$$(A31) \qquad \begin{aligned} &|(x, \hat{\Sigma}_{k,p} x) - ((J - \hat{E}_1)\hat{\Sigma}_{k,p}(J - \hat{E}_1)x, x)| \\ &\qquad \leq \|\hat{\Sigma}_{k,p}\|(2\|\hat{E}_1 - E_1\| + \|\hat{E}_1 - E_1\|^2). \end{aligned}$$

Therefore,

$$\begin{aligned} \lambda_2(\hat{\Sigma}_{k,p}) &= \max\{(x, (J - \hat{E}_1)\hat{\Sigma}_{k,p}(J - \hat{E}_1)x) : |x| = 1\} \\ &\leq O(\|\hat{E}_1 - E_1\|) + \lambda_2(\Sigma_p). \end{aligned}$$

Inverting the roles of $\hat{\Sigma}_{k,p}$ and $\Sigma_p$, we obtain

$$|\lambda_2(\hat{\Sigma}_{k,p}) - \lambda_2(\Sigma_p)| = O_P(\|\hat{\Sigma}_{k,p} - \Sigma_p\|).$$

Now repeating the argument we gave for (A29), we obtain

$$\|\hat{E}_2 - E_2\| = O_P(\|\hat{\Sigma}_{k,p} - \Sigma_p\|). \tag{A32}$$

The theorem follows from the induction and Theorem 1. $\square$

Note that if we track the effect of $\Delta$ and $m$, we in fact have

$$\|\hat{E}_j - E_j\| = O_P(j\Delta^{-2}\|\hat{\Sigma}_{k,p} - \Sigma_p\|), \qquad 1 \le j \le m.$$

Also note that the dimension of $\sum_{j=1}^m E_j$ is immaterial.

**Acknowledgments.** We thank Adam Rothman (University of Michigan) for help with simulations, Jianhua Huang (Texas A&M University) and Haipeng Shen (University of North Carolina at Chapel Hill) for sharing the call center data, and Sourav Chatterjee and Noureddine El Karoui (University of California, Berkeley) for helpful discussions. We also thank the two referees and the Co-editor, Jianqing Fan, for their helpful comments.

Department of Statistics  
University of California, Berkeley  
Berkeley, California 94720-3860  
USA  
E-mail: bickel@stat.berkeley.edu

Department of Statistics  
University of Michigan  
Ann Arbor, Michigan 48109-1107  
USA  
E-mail: elevina@umich.edu